\documentclass[aos,MSNbibl,nameyear,dvips]{arximspdf}
\usepackage{mathrsfs}
\usepackage{graphicx}

\doi{10.1214/13-AOS1149} 
\volume{41}
\issue{4}
\pubyear{2013}
\firstpage{2236}
\lastpage{2262}

\makeatletter

\newcommand{\rright}{\right}
\newcommand{\lleft}{\left}
\newcommand{\rrVert}{\Vert}

\newcommand{\llVert}{\Vert}

\newtheorem{theorem}{Theorem}
\newproclaim{definition}{Definition}
\newtheorem{lemma}{Lemma}
\newproclaim{condition}{Condition}

\newcommand{\bv}{\mathbf{v}}
\newcommand{\by}{\mathbf{y}}
\newcommand{\bA}{\mathbf{A}}
\newcommand{\bB}{\mathbf{B}}

\newcommand{\bQ}{\mathbf{Q}}
\newcommand{\bW}{\mathbf{W}}
\newcommand{\bX}{\mathbf{X}}
\newcommand{\bZ}{\mathbf{Z}}
\newcommand{\bzero}{\mathbf{0}}
\newcommand{\bbeta}{\bolds{\beta}}
\newcommand{\bSig}{\bolds{\Sigma}}

\newcommand{\Sig}{\bolds{\Sigma}}
\newcommand{\diag}{\operatorname{diag}}

\makeatother

\begin{document}
\begin{frontmatter}

\title{Impacts of high dimensionality in finite samples\thanksref{T1}}
\runtitle{Impacts of high dimensionality in finite samples}

\thankstext{T1}{Supported by NSF CAREER Award DMS-09-55316 and Grant
DMS-08-06030.}

\begin{aug}
\author[A]{\fnms{Jinchi} \snm{Lv}\corref{}\ead[label=e1]{jinchilv@marshall.usc.edu}}
\runauthor{J. Lv}
\affiliation{University of Southern California}
\address[A]{Data Sciences and Operations Department\\
Marshall School of Business\\
University of Southern California\\
Los Angeles, California 90089\\
USA\\
\printead{e1}} 
\end{aug}

\received{\smonth{6} \syear{2012}}
\revised{\smonth{3} \syear{2013}}

%
\begin{abstract}
High-dimensional data sets are commonly collected in many contemporary
applications arising in various fields of scientific research. We
present two views of finite samples in high dimensions: a probabilistic
one and a nonprobabilistic one. With the probabilistic view, we
establish the concentration property and robust spark bound for large
random design matrix generated from elliptical distributions, with the
former related to the sure screening property and the latter related to
sparse model identifiability. An interesting concentration phenomenon
in high dimensions is revealed. With the nonprobabilistic view, we
derive general bounds on dimensionality with some distance constraint
on sparse models. These results provide new insights into the impacts
of high dimensionality in finite samples.
\end{abstract}

%

\end{frontmatter}

\section{Introduction} \label{Sec1}
Thanks to the advances of information technologies, large-scale data
sets with a large number of variables or dimensions are commonly
collected in many contemporary applications that arise in different
fields of sciences, engineering and humanities. Examples include
marketing data in business, panel data in economics and finance,
genomics data in heath sciences, and brain imaging data in
neuroscience, among many others. The emergence of a large amount of
information contained in high-dimensional data sets provides
opportunities, as well as unprecedented challenges, for developing new
statistical methods and theory. See, for example, \citet{Hal06}
and \citet{FanLi06} for insights and discussions on the
statistical challenges associated with high dimensionality, and
\citet{FanLv10} for a brief review of some recent developments in
high-dimensional sparse modeling with variable selection. The approach
of variable selection aims to effectively identify important variables
and efficiently estimate their effects on a response variable of
interest.

For the purpose of prediction and variable selection, it is important
to understand and characterize the impacts of high dimensionality in
finite samples. \citet{HalMarNee05} investigated this
problem under the asymptotic framework of fixed sample size $n$ and
diverging dimensionality $p$, and revealed an interesting geometric
representation of high dimension, low sample size data. When viewed in
the diverging $p$-dimensional Euclidean space, the randomness in the
data vectors can be asymptotically squeezed into random rotation, with
the shape of the rescaled $n$-polyhedron approaching deterministic,
modulo the orientation. Such concentration phenomenon of random design
matrix in high dimensions is also shared by the concentration property
in \citet{FanLv08N1} (see Definition \ref{Def1}), in the asymptotic
setting of both $n$ and $p$ diverging. Geometrically, this property
means that the configuration of the $n$ sub-data vectors, modulo the
orientation, becomes close to regular asymptotically. Such a property
is key to establishing the sure screening property, which means that
all important variables are retained in the reduced feature space with
asymptotic probability one, of the sure independence screening (SIS)
method introduced in \citet{FanLv08N1}.

The SIS uses the idea of independence learning by applying
componentwise regression. Techniques of independence learning have been
widely used for variable ranking and screening. Recent work on variable
screening includes \citet{FanFan08}, \citet{HalTitXue09}, \citet{FanFenSon11},
\citet{XueZou11}, \citet{Zhuetal11}, \citet{DelHal12}, \citet{LiZhoZhu12},
\citet{MaiZou13} and \citet{BuhMan}, among others. The utility
of these methods is characterized by the sure screening property. In
particular, \citet{FanLv08N1} proved that the concentration property
holds when the
design matrix is generated from Gaussian distribution, and conjectured
that it may well hold for a wide class of elliptical distributions.
\citet{FanLv08N2} presented some simulation studies investigating such a
property for non-Gaussian distributions.
The first major contribution of our paper is to provide an affirmative
answer to the conjecture posed in \citet{FanLv08N1}.

To ensure model identifiability and stability for reliable
prediction and variable selection, it is practically important to control
the collinearity for sparse models. Since it is well known that the
level of collinearity among covariates typically increases with the
model dimensionality, bounding the sparse model size can be effective
in controlling model collinearity. Such a bound is characterized by the
concept of robust spark (see Definition \ref{Def2}). Another
contribution of the paper is to establish a lower bound on the robust
spark in the setting of large random design matrix generated from the
family of elliptical distributions.

In addition to the above probabilistic view of finite samples in high
dimensions, we also present a nonprobabilistic high-dimensional
geometric view. Both views are concerned with how much information
finite sample contains. A fundamental question is what the impact of
high dimensionality on differentiating the subspaces spanned by
different sets of predictors is. Such a question is tied to the issue
of model identifiability. In this paper, we intend
to derive general bounds on dimensionality with some distance
constraint on sparse models.

The rest of the paper is organized as follows. Section~\ref{Sec2}
establishes the concentration property and robust spark bound for large
random design matrix generated from elliptical distributions. We
investigate general bounds on dimensionality with distance constraint
from a nonprobabilistic point of view in Section~\ref{Sec3}. Section~\ref{Sec4} presents two numerical examples to illustrate the
theoretical results. We provide some discussions of our results and
their implications in Section~\ref{Sec5}. All technical details are
relegated to the \hyperref[app]{Appendix}.

\section{Concentration property and robust spark bound of large random
design matrix} \label{Sec2}
In this section, we focus on the case of large random design matrix
observed in a high-dimensional problem, in which each column vector
contains the information of a particular covariate. In high-dimensional
sparse modeling, a~common practice is to assume that only a faction of
all covariates, the so-called true or important covariates, contribute
to the regression or classification problem, whereas the other
covariates, the so-called noise covariates, are simply noise
information. The inclusion of noise covariates can deteriorate the
performance of the estimated model due to the well-known phenomenon of
noise accumulation in high-dimensional prediction
[\citet{FanFan08,FanLv10}]. A crucial issue behind high-dimensional inference is
to characterize the distance between the true underlying sparse model
and other sparse models, under some discrepancy measure. Intuitively,
such a distance can become smaller as the dimensionality increases,
making it more difficult to distinguish the true model from the others.
Therefore, it is a fundamental problem to characterize the impacts of
high dimensionality in finite samples.

\subsection{Concentration property} \label{Sec21}
We start with the task of dimensionality reduction, particularly
variable screening, which is useful in analyzing ultra-high dimensional
data sets. With the idea of independence learning, \citet{FanLv08N1}
introduced the SIS method to reduce the dimensionality of the feature
space from the ultra-high scale to a moderate scale, such as below
sample size. They introduced an asymptotic framework under which the
SIS enjoys the sure screening property even when the dimensionality can
grow exponentially with the sample size; see their Theorem 1. The sure
screening property means that the true model is contained in the much
reduced model after variable screening with asymptotic probability one.
In particular, a key ingredient of their asymptotic analysis is the
so-called concentration property in Condition 2 of \citet{FanLv08N1}.
They verified such property for random design matrix generated from
Gaussian distribution, and conjectured that it may also hold for a wide
class of elliptical distributions. To show that SIS is widely
applicable, it is crucial to establish the concentration property for
classes of non-Gaussian distributions.

The class of elliptical distributions, which is a wide family of
distributions generalizing the multivariate normal distribution, has
been broadly used in real applications. Examples of nonnormal
elliptical distributions include the Laplace distribution,
$t$-distribution, Cauchy distribution, logistic distribution, and
symmetric stable distribution. In particular, an important subclass of
elliptical distributions is the family of mixtures of normal
distributions. Mixture distributions provide a useful tool for
describing heterogeneous populations. Elliptical distributions also
play an important role in the theory of portfolio choice
[\citet{Cha83,OweRab83}]. This is due to important properties that any
affine transformation of elliptically distributed random variables
still has an elliptical distribution and each elliptical distribution
is uniquely determined by its location and scale parameters. An
implication in portfolio theory is that if all asset returns jointly
follow an elliptical distribution, then all portfolios are
characterized fully by their location and scale parameters. We refer to
\citet{FanKotNg90} for a comprehensive account of elliptical
distributions.

Assume that $\mathbf{x}= (X_1,\ldots, X_p)^T$ is a $p$-dimensional random
covariate vector having an elliptical distribution $\nu_p$ with mean
$\bzero$ and nonsingular covariance matrix~$\Sig$, and that we have a
sample $(\mathbf{x}_i)_{i = 1}^n$ of i.i.d. covariate vectors from this
distribution. Then we have an $n \times p$ random design matrix $\bX=
(\mathbf{x}_1,\ldots, \mathbf{x}_n)^T$. By the definition of
elliptical distribution
[see \citet{Mui82} or \citet{FanKotNg90}], the transformed
$p$-dimensional random vector $\mathbf{z}= (Z_1,\ldots, Z_p)^T=
\Sig^{-1/2}
\mathbf{x}$ has a spherical distribution $\mu_p$ with mean $\bzero
$ and
covariance matrix $I_p$. Similarly, we define the transformed covariate
vectors and transformed random design matrix as
%
%
\begin{equation}
\label{e001} \mathbf{z}_i = \Sig^{-1/2} \mathbf{x}_i
\quad\mbox{and}\quad \bZ= \bX\Sig^{-1/2},
\end{equation}
where $i = 1,\ldots, n$. Clearly, $\mathbf{z}_1,\ldots,
\mathbf{z}_n$ are $n$
i.i.d. copies of the transformed random covariate vector $\mathbf{z}$. We
denote by $\lambda_{\max}(\cdot)$ and $\lambda_{\min}(\cdot)$ the
largest and smallest eigenvalue of a given matrix, respectively. In
high-dimensional problems, we often face the situation of $p \gg n$, so
it is desirable to reduce the dimensionality of the feature space from
$p$ to a moderate one such as below sample size $n$. The SIS is capable
of doing so when the random design matrix $\bX$ satisfies the following
property, as introduced in \citet{FanLv08N1}.

\begin{definition}[(Concentration property)] \label{Def1}
The random design matrix $\bX$ is said to satisfy the concentration
property if there exist some positive constants $c_1, C_1$ such that
the deviation probability bound
%
%
\begin{equation}
\label{e002} P \bigl\{\lambda_{\max}\bigl(\widetilde{p}{}^{-1}
\widetilde{\bZ} \widetilde{\bZ}{}^T\bigr) > c_1  \mbox{
or }  \lambda_{\min}\bigl(\widetilde{p}{}^{-1} \widetilde{\bZ}
\widetilde{\bZ}{}^T\bigr) < c_1^{-1} \bigr\} \leq
\exp(-C_1 n)
\end{equation}
holds for each $n \times\widetilde{p}$ submatrix $\widetilde{\bZ}$ of
$\bZ$ with $c n < \widetilde{p} \leq p$ and $c > 1$ some positive constant.
\end{definition}

As mentioned in the \hyperref[Sec1]{Introduction}, the above
concentration property shows a similar concentration phenomenon of
large random design matrix to that in \citet{HalMarNee05}.
When the distribution $\nu _p$ of the covariate vector $\mathbf{x}$ is
$p$-variate Gaussian, \citet{FanLv08N1} proved that the random design
matrix $\bX$ satisfies the concentration property. We now consider a
more general class of distributions including Gaussian distributions,
the family of elliptical distributions. Assume that $P(\mathbf{z}=
\bzero) = 0$. Then it follows from Theorem 1.5.6 in \citet{Mui82}
that the $p$-variate spherical distribution $\mu_p$ has a density
function with respect to the Lebesgue measure that is spherically
symmetric on $\mathbb{R}^p$. We will work with the family of
log-concave spherical distributions on $\mathbb{R}^p$ that satisfy the
following two conditions.

\begin{condition} \label{cond1}
The density function $\exp\{-U_p(\bv)\}$ of the $p$-variate log-concave
spherical distribution $\mu_p$ satisfies that for some positive
constant~$c_2$,
%
%
\begin{equation}
\label{e003} \nabla^2 U_p(\bv) \geq c_2
I_p \qquad\mbox{uniformly in } \bv\in\mathbb{R}^p,
\end{equation}
where $\nabla^2 U_p$ denotes the Hessian matrix and $\bA\geq\bB$
means that $\bA- \bB$ is positive semidefinite for any symmetric
matrices $\bA$ and $\bB$.
\end{condition}

\begin{condition} \label{cond2}
There exists some positive constant $c_3 \leq1$ such that $E |Z_1|
\geq c_3$.
\end{condition}

Condition \ref{cond1} puts a constraint on the curvature of the
log-density of distribution~$\mu_p$, and Condition \ref{cond2} requires
that the mean $E |Z_1|$ needs to be bounded from below. Clearly,
log-concave spherical distributions satisfying Conditions \ref
{cond1}--\ref{cond2} comprise a wide class containing Gaussian
distributions. As seen in Lemma \ref{Lem2} later, Condition \ref{cond1}
entails that the corresponding spherical distribution is light-tailed,
which is important for variable screening. For heavy-tailed data sets,
\citet{DelHal12} showed that effective variable selection with
untransformed data requires slower growth of dimensionality. In
particular, they exploited variable transformation methods to transform
the original data into light-tailed data and demonstrated their
effectiveness and advantages. So in the presence of heavy-tailed data,
one may work with the assumption of elliptical distributions on the
transformed data.

The assumption of elliptical distributions is commonly used in
dimension reduction and has also been used for variable screening. See,
for example, \citet{Zhuetal11}. This assumption facilitates our
technical analysis. Similar results may hold for more general family of
distributions by resorting to techniques in random matrix theory. Some
other variable screening methods such as in \citet{MaiZou13} require
no such an assumption.

\begin{theorem}[(Concentration property)] \label{Thm1}
Under Conditions \ref{cond1}--\ref{cond2}, the random design matrix
$\bX
$ satisfies the concentration property (\ref{e002}).
\end{theorem}

Theorem \ref{Thm1} shows that the concentration property holds not only
for Gaussian distributions, but also for a wide class of elliptical
distributions, as conjectured by \citet{FanLv08N1} (see their Section~5.1). This provides an affirmative answer to their conjecture, showing
that the SIS indeed enjoys the sure screening property for the random
design matrix generated from a wide class of elliptical distributions.
The proof of Theorem \ref{Thm1} relies on the following three lemmas
that are of independent interest.

\begin{lemma} \label{Lem1}
Under Condition \ref{cond1}, each $q$-variate marginal distribution
$\widetilde{\mu}_q$ of $\mu_p$ with $1 \leq q \leq p$ satisfies the
logarithmic Sobolev inequality
%
%
\begin{equation}
\label{e004} E_{\widetilde{\mu}_q} \bigl\{f^2 \log f^2 \bigr
\} \leq2 C_2 E_{\widetilde{\mu}_q} \|\nabla f\|_2^2
\end{equation}
for any smooth function $f$ on $\mathbb{R}^q$ with $E_{\widetilde{\mu
}_q} f^2 = 1$, where $C_2 = c_2^{-1}$ and $\nabla f$ denotes the
gradient of function $f$.
\end{lemma}

\begin{lemma} \label{Lem2}
Let $\mathbf{z}_q$ be an arbitrary $q$-dimensional subvector of
$\mathbf{z}$ with \mbox{$1
\leq q \leq p$}. Then we have:
\begin{longlist}[(a)]
\item[(a)]
Under Condition \ref{cond1}, it holds for any $r \in(0, \infty)$ that
%
%
\begin{equation}
\label{e007} P \bigl\{\bigl|\|\mathbf{z}_q\|_2 - E
\|\mathbf{z}_q\|_2\bigr|> r \bigr\} \leq2 \exp\bigl
\{-C_2^{-1} r^2/2\bigr\};
\end{equation}

\item[(b)]
It holds that
%
%
\begin{equation}
\label{e008} \sqrt{q} E |Z_1| \leq E \|\mathbf{z}_q
\|_2 \leq\sqrt{q}.
\end{equation}
\end{longlist}
\end{lemma}

\begin{lemma} \label{Lem3}
Assume that Conditions \ref{cond1}--\ref{cond2} hold, $\mathbf{z}_q$ is a
$q$-dimensional subvector of $\mathbf{z}$ with $n \leq q \leq p$, and
$\mathbf{w}\sim
N(\bzero, I_q)$. Then there exist some positive constants $c_4 < 1$,
$c_5 > 1$, and $C_3$ such that
%
%
\begin{equation}
\label{e009} P \biggl\{\frac{\|\mathbf{z}_q\|_2}{\|\mathbf{w}\|_2}
< c_4  \mbox{ or }
\frac{\|\mathbf{z}
_q\|_2}{\|\mathbf{w}\|_2} > c_5 \biggr\} \leq4 \exp\{-C_3 n\}.
\end{equation}
\end{lemma}

Lemma \ref{Lem1} shows that each marginal distribution of $\mu_p$
satisfies the logarithmic Sobolev inequality, which is an important
tool for proving the concentration probability inequality for measures.
Lemma \ref{Lem2} establishes that for each $q$-dimensional subvector
$\mathbf{z}_q$ of $\mathbf{z}$, its $L_2$-norm $\|\mathbf{z}_q\|_2$ concentrates around the
mean $E \|\mathbf{z}_q\|_2$ with significant probability, which is
in turn
sandwiched between two quantities $\sqrt{q} E |Z_1|$ and $\sqrt{q}$.
Lemma \ref{Lem3} demonstrates an interesting phenomenon of measure
concentration in high dimensions.

\subsection{Robust spark bound} \label{Sec22}
As is well known in high-dimensional sparse modeling, controlling the
level of collinearity for sparse models is essential for model
identifiability and stable estimation. For a given $n \times p$ design
matrix $\bX$, there may exist another $p$-vector $\bbeta_1$ that is
different from the true regression coefficient vector $\bbeta_0$ such
that $\bX\bbeta_1$ is (nearly) identical to $\bX\bbeta_0$, when the
dimensionality $p$ is large compared with the sample size $n$. This
indicates that model identifiability is generally not guaranteed in
high dimensions when no additional constraint is imposed on the model
parameter. In addition, the subdesign matrix corresponding to a sparse
model should be well conditioned to ensure reliable estimation of model
parameters and nice convergence rates as in such as the least-squares
or maximum likelihood estimation. As an example, the covariance matrix
of the least-squares estimator is proportional to the inverse Gram
matrix given by the design matrix.

Since the collinearity among the covariates increases with the
dimensionality as evident from the geometric point of view, a natural
and effective way to ensure model identifiability and reduce model
instability is to control the size of sparse models. Such an idea has
been adopted in \citet{DonEla03} for the problem of sparse
recovery, which is the noiseless case of linear regression. In
particular, they introduced the concept of spark as a bound on sparse
model size to characterize model identifiability. The spark $\kappa$
of the
design matrix $\bX$ is defined as the smallest possible positive
integer such that there exists a singular $n \times\kappa$ submatrix
of $\bX$. This concept plays an important role in the problem of sparse
recovery; see also \citet{LvFan09}. An implication is that the true
model parameter vector $\bbeta_0$ is uniquely defined as long as $\|
\bbeta_0\|_0 < \kappa/2$, which provides a basic condition for model
identifiability. For the problem of variable selection in the presence
of noise, a stronger condition than provided by the spark is generally
needed. For this purpose, the concept of spark was generalized in
\citet{ZheFanLv} by introducing the concept of robust spark, as follows.

\begin{definition}[(Robust spark)] \label{Def2}
The robust spark $\kappa_c$ of the $n \times p$ design matrix $\bX$ is
defined as the smallest possible positive integer such that there
exists an $n \times\kappa_c$ submatrix of $n^{-1/2}\bX$ having a
singular value less than a given positive constant $c$.
\end{definition}

It is easy to see that the robust spark $\kappa_c$ approaches the spark
of $\bX$ as $c\rightarrow0+$. The robust spark provides a natural
bound on model size for effectively controlling the collinearity level
of sparse models, which is referred to as the robust spark condition.
For each sparse model with size $d < \kappa_c$, the corresponding $n
\times d$ submatrix of $n^{-1/2}\bX$ have all singular values bounded
from below by $c$. The robust spark $\kappa_c$ is always a positive
integer no larger than $n + 1$. It is practically important in
high-dimensional sparse modeling to show that the robust spark can be
some large number diverging with the sample size $n$. We intend to
build a lower bound on the robust spark for the case of random design
matrix, following the setting in Section~\ref{Sec21}.

\begin{theorem}[(Robust spark bound)] \label{Thm2}
Assume that the rows of the $n \times p$ random design matrix $\bX$ are
i.i.d. as $\nu_p$ having mean $\bzero$ and covariance matrix $\Sig$ and
satisfying Conditions \ref{cond1}--\ref{cond2}, with $\lambda_{\min
}(\bSig)$ bounded from below by some positive constant. Then with
asymptotic probability one, $\kappa_c \geq\widetilde{c} n/(\log p)$ for
sufficiently small constant $c$ and some positive constant $\widetilde{c}$
depending only on $c$.
\end{theorem}

Theorem \ref{Thm2} formally characterizes the order of the robust spark
$\kappa_c$ when the design matrix $\bX$ is generated from the family of
elliptical distributions. We see that sparse linear models of size as
large as of order $O \{n/(\log p) \}$ can still be well
separated from each other. On the other hand, when the true model size
is beyond such an order, the true underlying sparse model may be
indistinguishable from others in finite sample. Theorem \ref{Thm2} also
justifies the range of the true sparse model size under which the
problem of variable selection is meaningful. The deflation factor of
$\log p$ represents the general price one has to pay for the search of
important covariates in high dimensions.

The concept of robust spark shares a similar spirit as the restricted
eigenvalue condition on the design matrix in \citet{BicRitTsy09},
in the sense that both are sparse eigenvalue type conditions. Instead
of constraining the sparse model size, the restricted eigenvalue
condition uses an $L_1$-norm constraint on the parameter vector. As
discussed in \citet{ZheFanLv}, the robust spark condition can be
weaker than the restricted eigenvalue condition, since the $L_0$-norm
constraint can define a smaller subset than the $L_1$-norm constraint.
Many other conditions have also been introduced to characterize the
properties of variable selection methods such as the Lasso. See, for
example, \citet{vanBuh09} for a comprehensive
comparison and discussions on these conditions.

In particular, the robust spark condition is
weaker than the partial orthogonality condition, which requires that
true covariates and noise covariates are essentially uncorrelated, with
absolute correlation of the order $O(n^{-1/2})$. In contrast, the
robust spark condition can allow for much stronger correlation between
true covariates and noise covariates. The robust spark condition can
also be weaker than the irrepresentable condition. To see this, let us
consider the simple example constructed in \citet{LvFan09}. In their
Example 1, the irrepresentable condition becomes the constraint that
the maximum absolute correlation between the response and all noise
covariates is bounded from above by $s^{-1/2}$, where $s$ denotes the
true model size. Since the response is a linear combination of true
covariates in that example, this indicates that the irrepresentable
condition can be stronger than the robust spark condition when the true
model size grows.

\section{General bounds on dimensionality with distance constraint}
\label{Sec3}
We have provided in Section~\ref{Sec2} a probabilistic view of finite
samples in high dimensions, with focus on large random design matrix
generated from the family of elliptical distributions. It is also
important to understand how the dimensionality plays an role in
deterministic finite samples. For such a purpose, we take a
high-dimensional geometric view of finite samples and derive general
bounds on dimensionality using nonprobabilistic arguments. With a
slight abuse of notation, we now denote by $\mathbf{x}_j$ an $n$-dimensional
vector of observations from the $j$th covariate, and consider a
collection of $p$ covariates $\{\mathbf{x}_j\dvtx  j = 1,\ldots, p\}$.
Assume that
each covariate vector $\mathbf{x}_j$ is rescaled to have $L_2$-norm
$n^{1/2}$.
Then all vectors $n^{-1/2} \mathbf{x}_j$, $j = 1,\ldots, p$, lie
on the unit
sphere $S^{n - 1}$ in the $n$-dimensional Euclidean space $\mathbb
{R}^n$. We are interested in a natural question that how many variables
there can be if the maximum collinearity of sparse models is controlled.

For each positive integer $s$, denote by $\mathcal{A}_s$ the set of all
subspaces spanned by $s$ of covariates $\mathbf{x}_j$'s. Assume that
$s$ is
less than half of the spark $\kappa$ of the $n \times p$ design matrix
$\bX= (\mathbf{x}_1,\ldots, \mathbf{x}_p)$. Then each subspace
in $\mathcal{A}_s$
is $s$-dimensional and $|\mathcal{A}_s| = {p \choose s}$. To control
the collinearity among the variables, it is desirable to bound the
distances between $s$-dimensional subspaces in $\mathcal{A}_s$ away
from zero, under some discrepancy measure. When each pair of subspaces
in $\mathcal{A}_s$ has a positive distance, intuitively there cannot be
too many of them. The geometry of the space of all $s$-dimensional
subspaces of $\mathbb{R}^n$ is characterized by the Grassmann manifold
$G_{n, s}$. To facilitate our presentation, we list in Appendix \ref
{Sec31} some necessary background and terminology on the geometry and
invariant measure of Grassmann manifold. In particular, $G_{n, s}$
admits an invariant measure which under a change of variable and
symmetrization can be represented as a probability measure $\nu$ on
$[0, 1]^s$ with density given in (\ref{003}).

With the aid of the measure $\nu$, we can calculate the volumes of
various shapes of neighborhoods in the Grassmann manifold, which are
typically given in terms of the principal angles $\theta_i$ between an
$s$-dimensional subspace of $\mathbb{R}^n$ and a fixed $s$-dimensional
subspace with generator matrix $(I_s\enskip 0)$.
The principal angles between subspaces are natural extensions of the
concept of angle between lines. Let $V_1, V_2$ be two subspaces in
$G_{n, s}$ having a set of principal angles $(\theta_1,\ldots,
\theta
_s)$, with $\pi/2 \geq\theta_1 \geq\cdots\geq\theta_s \geq0$ and
corresponding $s$ pairs of unit vectors $(v_{1i}, v_{2i})$. If $V_i$ is
spanned by $s$ of $\mathbf{x}_j$'s, then putting $r_i = \cos\theta
_i$ and
reversing the order give the canonical correlations $(r_s,\ldots,
r_1)$ and corresponding pairs of canonical variables $(v_{1i},
v_{2i})$, for the two groups of variables.

There are three frequently used distances between subspaces $V_1$ and
$V_2$ on the Grassmann manifold $G_{n, s}$: the geodesic distance
$d_g(V_1, V_2) = (\sum_{i = 1}^s \theta_i^2)^{1/2}$ [\citet{Won67}], the
chordal distance $d_c(V_1, V_2) = (\sum_{i = 1}^s \sin^2 \theta
_i)^{1/2}$ [\citet{ConHarSlo96}], and the maximum chordal
distance [\citet{EdeAriSmi99}]
%
%
\begin{equation}
\label{004} d_m(V_1, V_2) = \sin
\theta_1 = \max_{i = 1}^s \sin
\theta_i.
\end{equation}
In view of the probability measure $\nu$ in (\ref{003}), it seems
natural to consider the latter two distances,
which is indeed the case. To see this, let $B_i$ be an $s \times n$
orthonormal generator matrix for $V_i$. Then $V_i$ is uniquely
determined by the projection matrix $P_i = B_i^T B_i$, which
corresponds to the projection onto the $s$-dimensional subspace~$V_i$.
It is known that
\[
d_c(V_1, V_2) = 2^{-1/2}
\|P_1 - P_2\|_F \quad\mbox{and}\quad
d_m(V_1, V_2) = \|P_1 -
P_2\|_2,
\]
where \mbox{$\|\cdot\|_F$} and \mbox{$\|\cdot\|_2$} denote the Frobenius norm and
spectral norm (or operator norm) of a given matrix, respectively. These
two matrix norms are commonly used in large covariance matrix
estimation and other multivariate analysis problems.

We now bound the size of the set $\mathcal{A}_s \subset G_{n, s}$ of
all subspaces spanned by $s$ of covariates $\mathbf{x}_j$'s under some
distance constraint, which in turn gives bounds on the dimensionality
$p$. The probability measure $\nu$ on $[0, 1]^s$ defined in (\ref{003})
is a key ingredient in our analysis. When all the subspaces in
$\mathcal
{A}_s$ have distance at least $2 \delta> 0$ under any distance $d$, it
is\vspace*{1pt} easy to see that ${p \choose s} = |\mathcal{A}_s| \leq1/\nu
(B_{\delta, d})$, where $B_{\delta, d}$ denotes a ball of radius
$\delta
$ in Grassmann manifold $G_{n, s}$ under distance $d$. In particular,
we focus on the maximum chordal distance defined in (\ref{004}).
Equivalently, the maximum chordal distance constraint gives the maximum
principal angle constraint. Since the sample size $n$ is usually small
or moderate in many contemporary applications, we adopt the asymptotic
framework of $s/n \rightarrow\gamma\in(0, 1)$ as $n \rightarrow
\infty$ for deriving asymptotic bounds on the dimensionality $p$.

\begin{theorem} \label{Thm3}
Assume that all subspaces spanned by $s$ of covariates $\mathbf{x}_j$'s have
maximum chordal distance
at least a fixed constant $2 \delta\in(0, 1)$, and $s/n \rightarrow
\gamma\in(0, 1/2)$ as $n \rightarrow\infty$. Then we have
%
%
\begin{equation}
\label{009} \log p \lesssim\bigl(\log\delta^{-1}\bigr) (1-\gamma)n + 2
\log n + O(1),
\end{equation}
where $\lesssim$ denotes asymptotic dominance.
\end{theorem}

Theorem \ref{Thm3} gives a general asymptotic bound on the
dimensionality $p$ under the maximum chordal distance
constraint, or equivalently, the maximum principal angle constraint. We
see that finite sample can allow for a large number of variables, in
which sparse models with size much smaller than sample size $n$ can
still be distinguishable from each other. The leading order in the
bound for $\log p$ is proportional to sample size $n$, with factors
$\log\delta^{-1}$ and $1-\gamma$. This result is reasonable because
larger $\delta$ means bigger separation of all $s$-dimensional
subspaces spanned by covariates $\mathbf{x}_j$'s, and large $\gamma
$ means
more such subspaces separated from each other, both cases leading to
tighter constraint on the growth of dimensionality $p$. It is
interesting that there are only two terms $O(\log n)$ and $O(1)$
following the leading order in the above bound on dimensionality.

The general bound on dimensionality with distance constraint in Theorem
\ref{Thm3} also shares some similarity with the lower bound $O \{
n/(\log p) \}$ on the robust spark in Theorem \ref{Thm2}, although
the former uses nonprobabilistic arguments with no distributional
assumption and the latter applies probabilistic arguments. The robust
spark provides a natural bound on sparse model size to control
collinearity for sparse models. Intuitively, when the dimensionality
$p$ grows with the sample size $n$, one expects tighter control on the
robust spark through a deflation factor of $\log p$. Similarly, the
upper bound on the logarithmic dimensionality $\log p$ in Theorem \ref
{Thm3} decreases with the minimum maximum chordal distance $2 \delta$
between sparse models through the factor $\log\delta^{-1}$. As
mentioned in Section~\ref{Sec22}, these sparse eigenvalue type
conditions play an important role in characterizing the variable
selection properties including the model selection consistency for
various regularization methods. Although the result in Theorem \ref
{Thm3} can be viewed as the bound for the worst case scenario, it
provides us caution and guidance on the growth of dimensionality in
real applications, particularly when variable selection is an important
goal in the studies.

In general, the robust spark $\kappa_c$ provides a stronger measure on
collinearity than the maximum chordal distance. To see this, assume
that $s < \kappa_c/2$ and let $V_1, V_2$ be two subspaces spanned by
two different sets of $s$ of covariates $\mathbf{x}_j$'s. Then the maximum
principal angle $\theta_1$ between $V_1$ and $V_2$ is the angle between
two vectors $v_i \in V_i$, where $v_i$ is a linear combination of the
corresponding set of covariate vectors $n^{-1/2} \mathbf{x}_j$
for each $i = 1, 2$. Since the union of these two sets of covariates
has cardinality bounded from above by $2 s < \kappa_c$, it follows from
the definition of the robust spark that the angle $\theta_1$ between
$v_1$ and $v_2$ is bounded from zero, which entails that the maximum
chordal distance between $V_1$ and $V_2$ is also bounded from zero.
Conversely, when two $s$-dimensional subspaces $V_1$ and $V_2$ has the
maximum chordal distance bounded from zero, the subdesign matrix
corresponding to covariates in the sets can still be singular.

We next consider a stronger distance constraint than in Theorem \ref
{Thm3}, where in addition, all disjoint subspaces in $\mathcal{A}_s$
have minimum principal angles at least $\arcsin\delta_1$ for some
$\delta_1 \in(0, \delta]$, with $\delta$ given in Theorem \ref{Thm3}.
Such disjoint subspaces are spanned by disjoint sets of $s$ of
covariates $\mathbf{x}_j$'s. In this case, it is natural to expect a tighter
bound on the dimensionality $p$.

\begin{theorem} \label{Thm4}
Assume that the conditions of Theorem \ref{Thm3} hold and all disjoint
subspaces have minimum principal angles at least a fixed constant
$\arcsin\delta_1$ with $\delta_1 \in(0, \delta]$. Then we have
%
%
\begin{eqnarray}
\label{014}
\log(p-s) &\lesssim&\bigl[\bigl(\log\delta^{-1}\bigr) (1-
\gamma) -c_{\delta
_1} - \gamma- 2^{-1}\log(1-2\gamma) \bigr] n\nonumber\\[-8pt]\\[-8pt]
&&{} + 2
\log n + O(1),\nonumber
\end{eqnarray}
where $c_{\delta_1} = 2^{-1} [\log(1-\delta_1^2)^{-1} ]
(1-\gamma) - 2^{-1}(1-\delta_1^2)^{-1} \delta_1^2(1-2\gamma)$.
\end{theorem}

Compared to the bound in Theorem \ref{Thm3}, Theorem \ref{Thm4} indeed
provides a tighter bound on the dimensionality $p$ due to the
additional distance constraint involving $\delta_1$. We are interested
in the asymptotic bound on the dimensionality when $\delta_1$ is near
zero. In this case, we have $c_{\delta_1} \sim\delta_1^2 \gamma/2$.
Observe that $\gamma+ 2^{-1}\log(1-2\gamma) \leq0$ and is of order
$O(\gamma^2)$. It is generally difficult to derive tight bounds over
the whole ranges of $\delta_1$ and $\gamma$. This is essentially due to
the challenge of obtaining a globally tight function bounding the
function $f_1$ defined in (\ref{015}) from above, while retaining
analytical tractability of evaluating the resulting integral.

We finally revisit the marginal correlation ranking, a widely used
technique for analyzing large-scale data sets, from a nonprobabilistic
point of view. Given a sample of size $n$, the maximum correlation of
noise covariates with the response variable can exceed the maximum
correlation of true covariates with the response variable when the
dimensionality $p$ is high. Here the correlation between two
$n$-vectors $\bv_1$ and $\bv_2$ is referred to as $\cos\theta= \bv
_1^T\bv_2/(\|\bv_1\|_2\|\bv_2\|_2)$, where $\theta$ is the angle
between them. It is important to understand the limit on the
dimensionality $p$ under which the above undesired phenomenon can happen.

\begin{theorem} \label{Thm5}
Let $r \in(0, 1)$ be the maximum absolute correlation between $s$ true
predictors $\mathbf{x}_j$ and response vector $\by$ in $\mathbb
{R}^n$ and
assume that all $p - s$ noise predictors $\mathbf{x}_j$ have absolute
correlations bounded by $\delta\in(0, 1)$. Then there exists a noise
predictor having absolute correlation with $\by$ larger than $r$ if
$\log(p - s) \geq2^{-1} \{\log[4/(1-\delta^2)] \} (n-1) +
2^{-1} \log n + O(1)$.
\end{theorem}

It is an interesting result that the above asymptotic bound on the
dimensionality $p$ depends only on $\delta\in(0, 1)$, and is
independent of the specific value of $r \in(0, 1)$. The condition on
the dimensionality is sufficient but not necessary in general, since
one can always add an additional noise predictor having absolute
correlation with $\by$ larger than $r$. Nevertheless, Theorem \ref
{Thm5} gives us a general limit on dimensionality even when one
believes that a majority of noise predictors have weak correlation with
the response variable.

Meanwhile, we also see from Theorem \ref{Thm5} that the dimensionality
$p$ generally needs to be large compared to the sample size $n$ such
that a noise predictor may have the highest correlation with the
response variable. This result is reflected in
a common feature of many variable selection procedures including
commonly used greedy algorithms, that is, initially selecting one
predictor with the highest correlation with the response variable. See,
for example, the LARS algorithm in \citet{Efretal04} and the LLA
algorithm in \citet{ZouLi08}. Such a variable, which gives a sparse
model with size one, commonly appears on the solution paths of many
regularization methods for high-dimensional variable selection.

\begin{figure}

\includegraphics{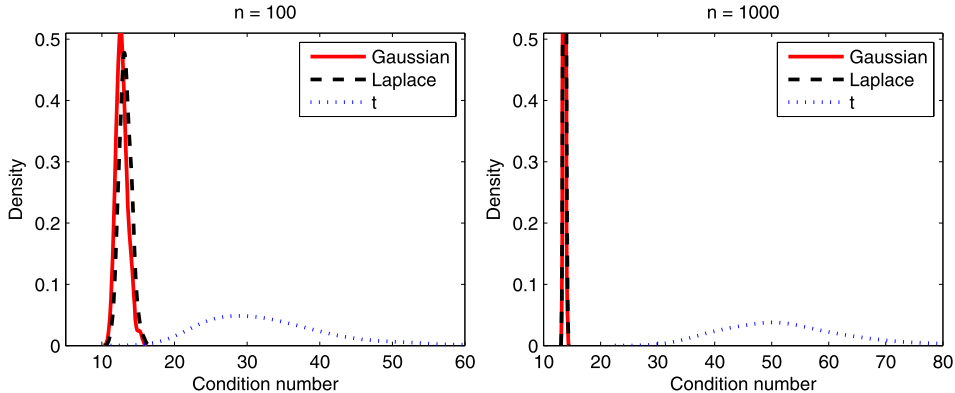}

\caption{Distributions of the condition number of $\widetilde{p}{}^{-1}
\bX\bX^T$ in different scenarios of distributions of $\bX$ with
$\widetilde{p} = 3 n$ and $n = 100$ and $1000$, based on 100
simulations.}\label{fig1}
\end{figure}

\section{Numerical examples} \label{Sec4}
In this section we provide two simulation examples to illustrate the
theoretical results in Section~\ref{Sec2}, obtained through
probabilistic arguments. The first simulation example examines the
concentration property for large random design matrix. Let $\bX$ be an
$n \times\widetilde{p}$ random design matrix with $\widetilde{p} = c
n$ for some constant $c > 1$. We set $n = 100$ and $1000$, and $c = 3$.
We considered three scenarios of distributions: (1) each entry of $\bX$
is sampled independently from $N(0, 1)$, (2) each entry of $\bX$ is
sampled independently from the Laplace distribution with mean 0 and
variance 1, and (3) each row of $\bX$ is sampled independently from the
multivariate $t$-distribution with 10 degrees of freedom and then
rescaled to have unit variances. In view of\vspace*{1pt} Definition \ref{Def1}, the
concentration property of $\bX$ is characterized by the distribution of
the condition number of $\widetilde{p}{}^{-1} \bX\bX^T$. In each case,
100 Monte Carlo simulations were used to obtain the distribution of
such condition number. Figure~\ref{fig1} depicts these distributions in
different scenarios. We see that in scenarios 1 and 2, the condition
number concentrates in the range of relatively small numbers,
indicating the associated concentration property as shown in Theorem
\ref{Thm1}. In scenario 3 with multivariate $t$-distribution, one still
observes the concentration phenomenon. However, since this distribution
is relatively more heavy-tailed, we see that the distribution of the
condition number becomes more spread out and shifts toward the range of
large numbers.

The second simulation example investigates the robust spark bound for
large \mbox{$n \times p$} random design matrix $\bX$. We adopted the same
three scenarios of distributions as in the first simulation example,
except that $n = 100$, and $p = 1000$ and $5000$. In light of Theorem
\ref{Thm2}, we sampled randomly 1000 $n \times k$ submatrices of
$n^{-1/2} \bX$ each with $k = \lceil2 n/(\log p)\rceil$ columns and
calculated the minimum of those 1000 smallest singular values.
Similarly, in each case 100 Monte Carlo simulations were used to obtain
the distribution of such minimum singular value which is tied to the
%
\begin{figure}

\includegraphics{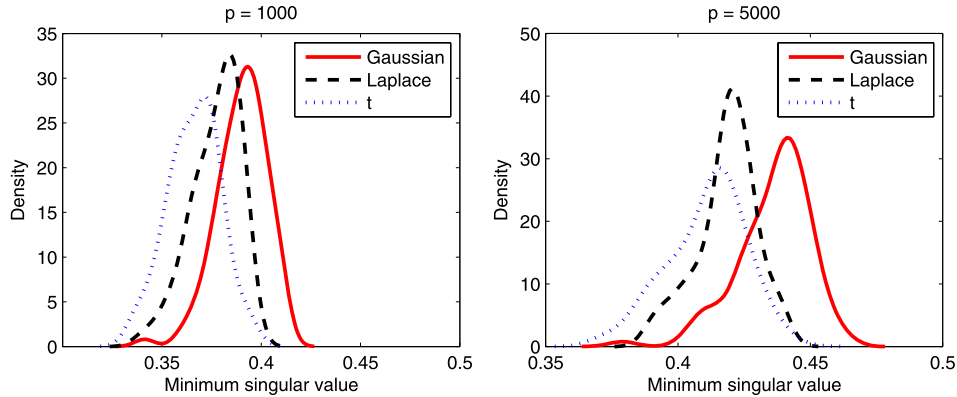}

\caption{Distributions of the minimum singular value over 1000
submatrices of $n^{-1/2} \bX$ each with $\lceil2 n/(\log p)\rceil$
columns, in different scenarios of distributions of $\bX$ with $n =
100$, and $p = 1000$ and $5000$, based on 100 simulations.}\label{fig2}
\end{figure}
robust spark bound of $\bX$. These distributions are shown in Figure~\ref{fig2}. In particular, we see that the distribution of the minimum
singular value concentrates clearly away from zero in each of the three
scenarios of distributions. These numerical results indicate that the
robust spark of random design matrix can indeed be at least of order
$O \{n/(\log p) \}$, as shown in Theorem \ref{Thm2}.

\section{Discussions} \label{Sec5}
We have investigated the impacts of high dimensionality in finite
samples from two different perspectives: a probabilistic one and a
nonprobabilistic one. An interesting concentration phenomenon for
large random design matrix has been revealed, as shown previously in
\citet{HalMarNee05}. We have shown that the concentration
property, which is important in characterizing the sure screening
property of the SIS, holds for a wide class of elliptical
distributions, as conjectured by \citet{FanLv08N1}. We have also
established a lower bound on the robust spark which is important in
ensuring model identifiability and stable estimation. The
high-dimensional geometric view of finite samples has lead to general
bounds on dimensionality with distance constraint on sparse models,
using nonprobabilistic arguments.

Both probabilistic and nonprobabilistic views provide understandings
on how the dimensionality interacts with the sample size for
large-scale data sets. Characterizing the limit of the dimensionality
with respect to the sample size is key to the success of
high-dimensional inference goals such as prediction and variable
selection. We have focused on the family of elliptical distributions.
It would be interesting to consider a more general class of
distributions for future research.

\begin{appendix}\label{app}
\section{Proofs of main results} \label{SecA}

For notational simplicity, we use $C$ to denote a generic positive
constant, whose value may change from line to line.

\subsection{\texorpdfstring{Proof of Lemma \protect\ref{Lem1}}{Proof of Lemma 1}} \label{SecA1}
By Theorem 5.2 in \citet{Led01}, we know that Condition \ref{cond1}
entails the logarithmic Sobolev inequality (\ref{e004}) when
$\widetilde
{\mu}_q = \mu_p$. It remains to prove the logarithmic Sobolev
inequality for any marginal distribution of $\mu_p$. Let $1 \leq q < p$
and $\widetilde{\mu}_q$ be a $q$-variate marginal distribution of
$\mu
_p$. By the spherical symmetry of $\mu_p$, without loss of generality
we can assume that $\widetilde{\mu}_q$ is concentrated on $\mathbb{R}^q
= \{\bv= (v_1,\ldots, v_p)^T\in\mathbb{R}^p\dvtx  v_{q + 1} = \cdots=
v_p = 0\}$. For any smooth function $\widetilde{f}$ on $\mathbb{R}^q$
with $E_{\widetilde{\mu}_q} \widetilde{f}{}^2 = 1$, define $f\dvtx  \mathbb
{R}^p \rightarrow\mathbb{R}$ by
%
%
\begin{equation}
\label{e005} f(v_1,\ldots, v_q, v_{q + 1},\ldots, v_p) = \widetilde{f}(v_1,\ldots,
v_q).
\end{equation}
Clearly $f$ is a smooth function on $\mathbb{R}^p$ and
\[
\nabla f(v_1,\ldots, v_q, v_{q + 1},\ldots,
v_p) = \lleft[ \matrix{ \nabla\widetilde{f}(v_1,\ldots, v_q)
\cr
\bzero} \rright],
\]
which shows that
%
%
\begin{equation}
\label{e006} \bigl\llVert\nabla f(v_1,\ldots, v_q,
v_{q + 1},\ldots, v_p)\bigr\rrVert_2^2
= \bigl\llVert\nabla\widetilde{f}(v_1,\ldots, v_q)
\bigr\rrVert_2^2.
\end{equation}
In view of (\ref{e005}), it follows from Fubini's theorem that
\[
E_{\mu_p} f^2 = E_{\widetilde{\mu}_q} \widetilde{f}{}^2
= 1.
\]
Thus by (\ref{e005}), (\ref{e006}), and Fubini's theorem, applying the
logarithmic Sobolev inequality (\ref{e004}) for $\mu_p$ to the smooth
function $f$ yields
\[
E_{\widetilde{\mu}_q} \bigl\{\widetilde{f}{}^2 \log\widetilde
{f}{}^2 \bigr\} = E_{\mu_p} \bigl\{f^2 \log
f^2 \bigr\} \leq2 C_2 E_{\mu_p} \| \nabla f
\|_2^2 = 2 C_2 E_{\widetilde{\mu}_q} \|\nabla
\widetilde{f}\|_2^2,
\]
which completes the proof.

\subsection{\texorpdfstring{Proof of Lemma \protect\ref{Lem2}}{Proof of Lemma 2}} \label{SecA2}
We first prove part (a). By Lemma \ref{Lem1}, the distribution of
$\mathbf{z}
_q$ satisfies the logarithmic Sobolev inequality (\ref{e004}). Observe
that by the triangle inequality, the Euclidean norm $\|\cdot\|_2$ is
$1$-Lipschitz with respect to the metric induced by itself. Therefore,
the classical Herbst argument applies to prove the concentration
inequality (\ref{e007}) [see, e.g., Theorem 5.3 in \citet{Led01}]. It
remains to show part (b). Note that $E Z_1^2 = 1$. By the spherical
symmetry of $\mu_p$, H\"{o}lder's inequality, and the Cauchy--Schwarz
inequality, we have
\[
E \|\mathbf{z}_q\|_2 \leq\sqrt{E \|\mathbf{z}_q\|_2^2} = \sqrt{q E
Z_1^2} = \sqrt{q}
\]
and
\[
E \|\mathbf{z}_q\|_2 \geq E \bigl(q^{-1/2} \|
\mathbf{z}_q\| _1 \bigr) = \sqrt{q} E
|Z_1|.
\]
This concludes the proof.

\subsection{\texorpdfstring{Proof of Lemma \protect\ref{Lem3}}{Proof of Lemma 3}} \label{SecA3}
We first make a simple observation. The standard Gaussian distributions
are special cases of spherical distributions. Recall that the
$q$-variate standard Gaussian distribution $\gamma_q$ has density
function\vspace*{1pt} $\frac{1}{(2 \pi)^{q/2}} e^{- \|\bv\|_2^2/2}$, $\bv\in
\mathbb
{R}^q$. Thus it\vspace*{1pt} is easy to check that $\gamma_q$ satisfies Condition
\ref{cond1} with $c_2 = 1$. Let $\mathbf{w}= (W_1,\ldots,
W_q)^T\sim
N(\bzero, I_q)$. Then it follows immediately from Lemma \ref{Lem2} that
for any $r \in(0, \infty)$,
%
%
\begin{equation}
\label{e010} P \bigl(\bigl|\|\mathbf{w}\|_2 - E \|\mathbf{w}
\|_2\bigr|> r \bigr) \leq2 e^{-r^2/2}.
\end{equation}
Note that
%
%
\begin{equation}
\label{e011} E |W_1| = 2 \int_0^\infty
u \frac{1}{\sqrt{2 \pi}} e^{-u^2/2} \,du = \sqrt{\frac{2}{\pi}} \int
_0^\infty e^{-u^2/2} \,d \biggl(
\frac
{u^2}{2} \biggr) = \sqrt{\frac{2}{\pi}}.
\end{equation}
By (\ref{e008}), (\ref{e010}), and (\ref{e011}), we have for any $r_1
\in(0, \sqrt{\frac{2}{\pi}})$,
%
%
\begin{eqnarray}
\label{e012} P \biggl(q^{-1/2} \|\mathbf{w}\|_2 > 1 +
r_1 \mbox{ or } q^{-1/2} \|\mathbf{w}\|_2 <
\sqrt{\frac{2}{\pi}} - r_1 \biggr) &\leq&2 e^{-q r_1^2/2}\nonumber\\[-8pt]\\[-8pt]
&\leq&2
e^{-n r_1^2/2}\nonumber
\end{eqnarray}
since $q \geq n$.

Now we get back to $\mathbf{z}_q$. It follows from (\ref{e007}) and
(\ref
{e008}) in Lemma \ref{Lem2} and Condition~\ref{cond2} that for any $r_2
\in(0, c_3)$,
%
%
\begin{eqnarray}
\label{e013} P \bigl(q^{-1/2} \|\mathbf{z}_q
\|_2 > 1 + r_2  \mbox{ or } q^{-1/2} \|\mathbf{z}_q\| _2 < c_3 - r_2 \bigr)
&\leq&2 e^{-C_2^{-1} q r_2^2/2} \nonumber\\[-8pt]\\[-8pt]
&\leq&2 e^{-C_2^{-1} n r_2^2/2}\nonumber
\end{eqnarray}
since $q \geq n$. Let
\[
c_4 = \frac{c_3 - r_2}{1 + r_1} \quad\mbox{and}\quad c_5 =
\frac{1
+ r_2}{\sqrt{{2}/{\pi}} - r_1}.
\]
Then combining (\ref{e012}) and (\ref{e013}) along with Bonferroni's
inequality yields
\begin{eqnarray*}
P \biggl(\frac{\|\mathbf{z}_q\|_2}{\|\mathbf{w}\|_2} < c_4  \mbox{
or }  \frac{\|\mathbf{z}
_q\|_2}{\|\mathbf{w}\|_2}
> c_5 \biggr) & \leq & 2 e^{-n r_1^2/2} + 2 e^{-C_2^{-1} n r_2^2/2}
\\
& \leq & 4 e^{-C_3 n},
\end{eqnarray*}
where $C_3 = \min(r_1^2/2, C_2^{-1} r_2^2/2)$. This completes the proof.

\subsection{\texorpdfstring{Proof of Theorem \protect\ref{Thm1}}{Proof of Theorem 1}} \label{SecA4}
In Section A.7, \citet{FanLv08N1} proved that Gaussian distributions
satisfy the concentration property (\ref{e002}), that is, for $\bZ
\sim
N(\bzero, I_n \otimes I_p) = N(\bzero, I_{n \times p})$. We now
consider the general situation where $n$ rows of the $n \times p$
random matrix $\bZ= (\mathbf{z}_1,\ldots, \mathbf{z}_n)^T$ are
i.i.d. copies from
the spherical distribution $\mu_p$. Fix an arbitrary $n \times
\widetilde{p}$ submatrix $\widetilde{\bZ}$ of $\bZ$ with $c n <
\widetilde{p} \leq p$, where $c \in(1, \infty)$. We aim to prove
deviation inequality in (\ref{e002}) with different constants $c_1 \in
(1, \infty)$ and $C_1 \in(0, \infty)$.

By the spherical symmetry, without loss of generality we can assume
that $\widetilde{\bZ}$ consists of the first $\widetilde{p}$ columns of
$\bZ$. Let
\[
\widetilde{\mathbf{z}} = (Z_1,\ldots, Z_{\widetilde{p}})^T
\quad\mbox{and}\quad \widetilde{\bZ} = (\widetilde{\mathbf{z}}_1,\ldots,
\widetilde{\mathbf{z}}_n)^T.
\]
Clearly, $\widetilde{\mathbf{z}}_1,\ldots, \widetilde{\mathbf{z}}_n$ are $n$ i.i.d.
copies of $\widetilde{\mathbf{z}}$. Take an $n \times\widetilde
{p}$ random matrix
\[
\bW= (\mathbf{w}_1,\ldots, \mathbf{w}_n)^T
\sim N(\bzero, I_n \otimes I_{\widetilde
{p}}),
\]
which is independent of $\widetilde{\bZ}$. Then for each
$i = 1,\ldots, n$, $\mathbf{w}_i$ has the $\widetilde{p}$-variate standard Gaussian
distribution. It is well known that $\mathbf{w}_i/\|\mathbf{w}_i\|
_2$ has the Haar
distribution on the unit sphere $S^{\widetilde{p} - 1}$ in $\widetilde
{p}$-dimensional Euclidean space $\mathbb{R}^{\widetilde{p}}$, that is,
the uniform distribution on $S^{\widetilde{p} - 1}$.

Since the distribution of $\widetilde{\mathbf{z}}$ is a marginal
distribution
of $\mu_p$, the spherical symmetry of $\mu_p$ entails that of the
distribution of $\widetilde{\mathbf{z}}$. It follows easily from the
assumption of $P(\mathbf{z}= \bzero) = 0$ that $P(\widetilde{
\mathbf{z}} = \bzero) =
0$. Thus by Theorem 1.5.6 in \citet{Mui82}, $\widetilde{\mathbf{z}}/\|
\widetilde{\mathbf{z}}\|_2$ is uniformly distributed on
$S^{\widetilde{p} -
1}$ and is independent of $\|\widetilde{\mathbf{z}}\|_2$. This
along with the
above fact shows that for each $i = 1,\ldots, n$,
%
%
\begin{equation}
\label{e014} \widetilde{\mathbf{z}}_i \stackrel{\mathrm{(d)}} {=}
\|\widetilde{\mathbf{z}}_i\|_2 \biggl(
\frac{\mathbf{w}_i}{\|\mathbf{w}
_i\|_2} \biggr) = \biggl(\frac{\|\widetilde{\mathbf{z}}_i\|_2}{\|
\mathbf{w}_i\|_2} \biggr) \mathbf{w}_i,
\end{equation}
where we use the symbol $\stackrel{\mathrm{(d)}}{=}$ to denote being
identical in
distribution. Hereafter, for notational simplicity we do not
distinguish $\widetilde{\mathbf{z}}_i$ and $ (\frac{\|
\widetilde{\mathbf{z}}_i\|
_2}{\|\mathbf{w}_i\|_2} ) \mathbf{w}_i$.

Define the $n \times n$ diagonal matrix
\[
\bQ= \diag\biggl\{\frac{\|\widetilde{\mathbf{z}}_1\|_2}{\|
\mathbf{w}_1\|_2},\ldots, \frac{\|\widetilde{\mathbf{z}}_n\|_2}{\|
\mathbf{w}_n\|_2} \biggr\}.
\]
Then we have
\[
\widetilde{p}{}^{-1} \widetilde{\bZ} \widetilde{\bZ}{}^T= \bQ
\bigl(\widetilde{p}{}^{-1} \bW\bW^T \bigr) \bQ,
\]
which entails
\begin{eqnarray*}
\min_{i = 1}^n \frac{\|\widetilde{\mathbf{z}}_i\|
_2^2}{\|\mathbf{w}_i\|_2^2}
\lambda_{\mathrm{min}}\bigl(\widetilde{p}{}^{-1} \bW\bW^T
\bigr) I_n & \leq & \min_{i = 1}^n
\frac{\|\widetilde{\mathbf{z}}_i\|_2^2}{\|
\mathbf{w}_i\|_2^2} \bigl(\widetilde{p}{}^{-1} \bW\bW^T \bigr)
\leq\bQ\bigl(\widetilde{p}{}^{-1} \bW\bW^T \bigr) \bQ
\\
& = & \widetilde{p}{}^{-1} \widetilde{\bZ} \widetilde{\bZ}{}^T
\leq\max_{i = 1}^n \frac{\|\widetilde{\mathbf{z}}_i\|_2^2}{\|
\mathbf{w}_i\|_2^2} \bigl(
\widetilde{p}{}^{-1} \bW\bW^T \bigr)
\\
& \leq &\max_{i = 1}^n \frac{\|\widetilde{\mathbf{z}}_i\|
_2^2}{\|\mathbf{w}
_i\|_2^2}
\lambda_{\mathrm{max}}\bigl(\widetilde{p}{}^{-1} \bW\bW^T
\bigr) I_n.
\end{eqnarray*}
This shows that
%
%
\begin{equation}
\label{e015} \lambda_{\mathrm{min}}\bigl(\widetilde{p}{}^{-1}
\widetilde{\bZ} \widetilde{\bZ}{}^T\bigr) \geq\min_{i = 1}^n
\frac{\|\widetilde{\mathbf{z}}_i\|_2^2}{\|\mathbf{w}
_i\|_2^2} \lambda_{\mathrm{min}}\bigl(\widetilde{p}{}^{-1} \bW
\bW^T\bigr)
\end{equation}
and
%
%
\begin{equation}
\label{e016} \lambda_{\mathrm{max}}\bigl(\widetilde{p}{}^{-1}
\widetilde{\bZ} \widetilde{\bZ}{}^T\bigr) \leq\max_{i = 1}^n
\frac{\|\widetilde{\mathbf{z}}_i\|_2^2}{\|\mathbf{w}
_i\|_2^2} \lambda_{\mathrm{max}}\bigl(\widetilde{p}{}^{-1} \bW
\bW^T\bigr).
\end{equation}

As mentioned before, we have for some $c_1 \in(1, \infty)$ and $C_1
\in(0, \infty)$,
%
%
\begin{equation}
\label{e017} P \bigl(\lambda_{\mathrm{max}}\bigl(\widetilde{p}{}^{-1}
\bW\bW^T\bigr) > c_1 \mbox{ or } \lambda_{\mathrm{min}}
\bigl(\widetilde{p}{}^{-1} \bW\bW^T\bigr) < 1/c_1
\bigr) \leq e^{-C_1n}.
\end{equation}
Note that $\widetilde{p} > n$. Thus by (\ref{e009}) in Lemma \ref
{Lem3}, an application of Bonferroni's inequality gives
%
%
\begin{equation}
\label{e018} P \biggl(\min_{i = 1}^n
\frac{\|\widetilde{\mathbf{z}}_i\|_2^2}{\|\mathbf{w}
_i\|_2^2} < c_4  \mbox{ or }  \max_{i = 1}^n
\frac{\|
\widetilde{\mathbf{z}}_i\|_2^2}{\|\mathbf{w}_i\|_2^2} > c_5 \biggr)
\leq4 n e^{-C_3 n},
\end{equation}
where $c_4 \in(0, 1)$, $c_5 \in(1, \infty)$, and $C_3 \in(0, \infty
)$. Therefore by Bonferroni's inequality, combining (\ref{e017}) and
(\ref{e018}) proves the deviation inequality in (\ref{e002}) by
appropriately changing the constants $c_1 \in(1, \infty)$ and $C_1
\in
(0, \infty)$. This concludes the proof.

\subsection{\texorpdfstring{Proof of Theorem \protect\ref{Thm2}}{Proof of Theorem 2}} \label{SecA5}
Using the similar arguments as in the proof of Theorem \ref{Thm1}, we
can prove that there exist some universal positive constants $c_6, C$
such that the deviation probability bound
%
%
\begin{equation}
\label{028} P \bigl\{\lambda_{\min}\bigl(n^{-1} \widetilde{
\bZ}{}^T\widetilde{\bZ}\bigr) < c_6 \bigr\} \leq\exp(-C n)
\end{equation}
holds for each $n \times\widetilde{p}$ submatrix $\widetilde{\bZ}$ of
$\bZ$ with $\widetilde{p} = c_7 n$ and $c_7 < 1$ some positive
constant. This is because Lemmas \ref{Lem1}--\ref{Lem2} are free of the
dimension $q$ of the marginal distribution, and Lemma \ref{Lem3} still
holds with the choice of $q = \widetilde{p}$. We should also note that
the deviation probability bound (\ref{028}) holds when $\widetilde
{\bZ}
\sim N(\bzero, I_n \otimes I_{\widetilde{p}})$, which is entailed by
the concentration property (\ref{e002}) proved in \citet{FanLv08N1} for
Gaussian distributions.

For each set $\alpha\subset\{1,\ldots, p\}$ with $|\alpha| =
\widetilde{p}$, denote by $\bSig_{\alpha,\alpha}$ the principal
submatrix of $\bSig$ corresponding to variables in $\alpha$, and $\bX
_\alpha$ a submatrix of the design matrix $\bX$ consisting of columns
with indices in $\alpha$. It follows easily from the representation of
elliptical distributions that $\widetilde\bX_{\alpha} =\bX_{\alpha
}\bSig
_{\alpha,\alpha}^{-1/2}$ has the same distribution as $\bZ_\alpha$.
Since $\lambda_{\min}(\bSig)$ is bounded from below by some positive
constant, we have
\[
\lambda_{\min} \bigl(n^{-1}\bX_{\alpha}^T
\bX_{\alpha} \bigr) \geq\lambda_{\min} \bigl(n^{-1}
\widetilde\bX_{\alpha}^T\widetilde\bX_{\alpha
} \bigr)
\lambda_{\min}(\bSig_{\alpha,\alpha})\geq C\lambda_{\min}
\bigl(n^{-1}\widetilde\bX_{\alpha}^T\widetilde
\bX_{\alpha} \bigr),
\]
where $C$ is some positive constant. Therefore, combining the above
results yields
%
%
\begin{equation}
\label{029} P \bigl\{\lambda_{\min} \bigl(n^{-1}
\bX_{\alpha}^T\bX_{\alpha} \bigr) < c_6
\bigr\} \leq\exp(-C n)
\end{equation}
with a possibly different positive constant $c_6$. Note that the
positive constants involved are universal ones. We choose a positive
integer $K = 2^{-1} C n/\break(\log p) \leq\widetilde{p}$. Then an
application of the Bonferroni inequality together with (\ref{029}) gives
\[
P \Bigl\{\min_{|\alpha| = K} \lambda_{\min}
\bigl(n^{-1}\bX_{\alpha
}^T\bX_{\alpha} \bigr)
< c_6 \Bigr\} \leq\sum_{|\alpha| = K}\exp(-C n)
\leq p^{K}\exp(-C n) \rightarrow0
\]
as $n \rightarrow\infty$. This shows that with asymptotic probability
one, the robust spark $\kappa_c \geq K$ for any $c \leq c_6$, which
completes the proof.

\subsection{\texorpdfstring{Proof of Theorem \protect\ref{Thm3}}{Proof of Theorem 3}} \label{SecA7}
For the maximum chordal distance $d_m$, by noting that $x_i = \sin^2
\theta_i$, we have a simple representation of the neighborhood
$B_{\delta, d_m} = \{(x_1,\ldots, x_s) \in[0, 1]^s\dvtx  \max_{i = 1}^s
x_i \leq\delta^2\}$ for $\delta\in(0, 1/2)$. We need to calculate
its volume under the probability measure $\nu$ given in (\ref{003}). In
light of (\ref{003}), a change of variable $y_i = \delta^{-2} x_i$ gives
%
%
\begin{equation}
\label{005}\quad d \nu= \delta^{s (n - s)} K_{n, s} f(y_1,\ldots, y_s) \prod_{1
\leq i
< j \leq s}
|y_i - y_j| \prod_{i = 1}^s
y_i^{\alpha-1} \,d y_1 \cdots d y_s,
\end{equation}
where $f(y_1,\ldots, y_s) = \prod_{i = 1}^s (1 - \delta^2 y_i)^{-1/2}$
over $[0, 1]^s$. Observe that without the term $f$ in (\ref{005}),
$\nu
(B_{\delta, d_m})$ would become Selberg's integral which is a
generalization of the beta integral [\citet{Meh04}]. We will evaluate
this integral by sandwiching the function $f$ between two functions of
the same form. Since the function $(1 - \delta^2 y)^{-1/2}$ is
increasing and convex on $[0, 1]$, it follows that $1 + c_1 y \leq(1 -
\delta^2 y)^{-1/2} \leq1 + c_2 y$, where $c_1 = \delta^2/2$ and $c_2 =
\delta^2 (1 - \delta^2)^{-3/2}/2$. This shows that
\[
\prod_{i = 1}^s (1 + c_1
y_i) \leq f(y_1,\ldots, y_s) \leq\prod
_{i =
1}^s (1 + c_2
y_i).
\]
Thus, we obtain a useful representation of the volume of the neighborhood
%
%
\begin{equation}
\label{008} \nu(B_{\delta, d_m}) = \delta^{s (n - s)} K_{n, s}
I(c),
\end{equation}
where $I(c)$ with some $c \in[c_1, c_2]$ is an integral given in the
following lemma.

\begin{lemma} \label{Lem4}
For each $c > 0$, we have
%
%
\begin{eqnarray}
\label{007} I(c) & \equiv & \int_{[0, 1]^s} \prod
_{i = 1}^s (1 + c y_i) \prod
_{1
\leq
i < j \leq s} |y_i - y_j| \prod
_{i = 1}^s y_i^{\alpha-1} \,d
y_1 \cdots d y_s
\nonumber
\\
& = & 2^s \pi^{-s/2} \sum_{m = 0}^s
\pmatrix{s \cr m} c^m \prod_{i = s -
m}^{s-1}
\frac{\alpha+ 2^{-1}i}{\alpha+ 2^{-1} (s + i + 1)}
\\
&&{} \times\prod_{i = 0}^{s - 1}
\frac{\Gamma(\alpha+ 2^{-1}i)
\Gamma(1 + 2^{-1} (i + 1)) \Gamma(1 + 2^{-1}i)}{\Gamma(\alpha+
2^{-1}(s + i +1))},
\nonumber
\end{eqnarray}
where the factor containing $m$ equals $1$ when $m = 0$.
\end{lemma}

\begin{pf}
Observe that the integrand in (\ref{007}) is symmetric in $y_1,\ldots,
y_s$. Thus, an expansion of $\prod_{i = 1}^s (1 + c y_i)$ gives
\begin{eqnarray*}
&& \int_{[0, 1]^s} \prod_{i = 1}^s
(1 + c y_i) \prod_{1 \leq i < j
\leq
s}
|y_i - y_j| \prod_{i = 1}^s
y_i^{\alpha-1} \,d y_1 \cdots d y_s
\\
&&\qquad = \sum_{m = 0}^s \pmatrix{s \cr m}
c^m \int_{[0, 1]^s} \prod
_{i = 1}^m y_i \prod
_{1 \leq i < j \leq s} |y_i - y_j| \prod
_{i = 1}^s y_i^{\alpha
-1} \,d
y_1 \cdots d y_s,
\end{eqnarray*}
where $\prod_{i = 1}^m y_i = 1$ when $m = 0$. The above integrals are
exactly Aomoto's extension of Selberg's integral [\citet{Meh04}] and can
be calculated as
\begin{eqnarray*}
&& \int_{[0, 1]^s} \prod_{i = 1}^m
y_i \prod_{1 \leq i < j \leq s} |y_i -
y_j| \prod_{i = 1}^s
y_i^{\alpha-1} \,d y_1 \cdots d y_s
\\
&&\qquad = 2^s \pi^{-s/2} \prod_{i = s - m}^{s-1}
\frac{\alpha+
2^{-1}i}{\alpha+ 2^{-1} (s + i + 1)} \\
&&\qquad\quad{}\times\prod_{i = 0}^{s - 1}
\frac
{\Gamma
(\alpha+ 2^{-1}i) \Gamma(1 + 2^{-1} (i + 1)) \Gamma(1 +
2^{-1}i)}{\Gamma(\alpha+ 2^{-1}(s + i +1))},
\end{eqnarray*}
where the factor containing $m$ equals $1$ when $m = 0$. This completes
the proof of Lemma \ref{Lem4}.
\end{pf}

Let us continue with the proof of Theorem \ref{Thm3}. By assumption, $s
\sim\gamma n$ as $n \rightarrow\infty$, so $\delta^{s (n - s)} \pi
^{-s/2} \sim\delta^{\gamma(1-\gamma)n^2} \pi^{-\gamma n/2}$. Applying
Stirling's formula for large factorials gives $s! \sim(2 \pi\gamma
n)^{1/2} (\gamma/e)^{\gamma n} n^{\gamma n}$. Thus by omitting $O(n)$
and smaller order terms,
%
%
\begin{equation}
\label{010} \log\bigl[\delta^{s (n - s)} \pi^{-s/2}/s!\bigr] \sim(
\log\delta) \gamma(1-\gamma)n^2 - \gamma n \log n.
\end{equation}
Using Stirling's formula for the Gamma function $\Gamma(t + 1) \sim(2
\pi t)^{1/2} (t/e)^t$ as $t \rightarrow\infty$ and noting that $A_j =
2 \pi^{j/2}/\Gamma(j/2)$ and $\alpha= (n - 2 s + 1)/2$, we derive
\begin{eqnarray*}
&&
\prod_{i = 0}^{s - 1} \frac{A_{n - s - i}}{A_{n - i}} \prod
_{i = 0}^{s
- 1} \frac{\Gamma(\alpha+ 2^{-1}i)}{\Gamma(\alpha+ 2^{-1}(s + i +1))}\\
&&\qquad
\sim
\pi^{-s^2/2} (2 e)^{s/2} (n-s-1)^{(n-s-1)/2} (n -
1)^{-(n - 1)/2},
\end{eqnarray*}
which entails that
%
%
\begin{eqnarray}
\label{011}
&&\log\Biggl\{\prod_{i = 0}^{s - 1}
\frac{A_{n - s - i}}{A_{n - i}} \prod_{i = 0}^{s - 1}
\frac{\Gamma(\alpha+ 2^{-1}i)}{\Gamma(\alpha+
2^{-1}(s + i +1))} \Biggr\} \nonumber\\[-8pt]\\[-8pt]
&&\qquad\sim-(\log\pi) \gamma^2 n^2/2
- \gamma n (\log n)/2.\nonumber
\end{eqnarray}
Similarly, it follows from the identities $\Gamma(t + 1)=t\Gamma(t)$
and $\Gamma(1)=1$ that
\[
2^{-s} \prod_{i = 0}^{s - 1}
A_{s - i}^2 \prod_{i = 0}^{s - 1}
\Gamma\bigl(1 + 2^{-1} (i + 1)\bigr) \Gamma\bigl(1 + 2^{-1}i
\bigr) = s! \pi^{s(s+1)/2}/\Gamma(s/2).
\]
This shows that
%
%
\begin{eqnarray}
\label{006}
&&\log\Biggl\{2^{-s} \prod_{i = 0}^{s - 1}
A_{s - i}^2 \prod_{i =
0}^{s -
1}
\Gamma\bigl(1 + 2^{-1} (i + 1)\bigr) \Gamma\bigl(1 + 2^{-1}i
\bigr) \Biggr\} \nonumber\\[-8pt]\\[-8pt]
&&\qquad\sim(\log\pi) \gamma^2 n^2/2 + \gamma n
(\log n)/2.\nonumber
\end{eqnarray}

It remains to consider the last term. Note that
\begin{eqnarray*}
&& \prod_{i = s - m}^{s-1} \frac{\alpha+ 2^{-1}i}{\alpha+ 2^{-1} (s +
i + 1)} \\
&&\qquad=
\frac{(n-s)!}{(n+1)!} \frac{(n-m+1)!}{(n-m-s)!}
\\
&&\qquad \sim\biggl(\frac{n-m+1}{n+1} \biggr)^{s+1} \biggl(
\frac
{n-s}{n+1} \biggr)^{n-s+1/2} \biggl(\frac{n-m+1}{n-m-s}
\biggr)^{n-m-s+1/2} \\
&&\qquad\sim e^{O(n)},
\end{eqnarray*}
which entails that
%
%
\begin{eqnarray}
\label{012} \log\Biggl\{\sum_{m = 0}^s \pmatrix{s
\cr m} c^m \prod_{i = s - m}^{s-1}
\frac{\alpha+ 2^{-1}i}{\alpha+ 2^{-1} (s + i + 1)} \Biggr\} &\sim&\log
\bigl[e^{O(n)} (1 +
c)^s\bigr] \nonumber\\[-8pt]\\[-8pt]
&=& O(n).\nonumber
\end{eqnarray}
Thus combining (\ref{008}) and (\ref{010})--(\ref{012}) yields
%
%
\begin{equation}
\label{013} \log\bigl(\nu(B_{\delta, d_m})\bigr) \sim(\log\delta) \gamma(1-
\gamma)n^2 - \gamma n \log n + O(n).
\end{equation}
Since all the subspaces in $\mathcal{A}_s$ have maximum chordal distance
at least $2 \delta$, it holds that ${p \choose s} = |\mathcal{A}_s|
\leq1/\nu(B_{\delta, d_m})$. Hence by (\ref{013}),
\[
\log\pmatrix{p \cr s} \lesssim\bigl(\log\delta^{-1}\bigr) \gamma(1-
\gamma)n^2 + \gamma n \log n + O(n),
\]
where $\lesssim$ denotes asymptotic dominance. It is easy to derive
$\log{p \choose s} \gtrsim\gamma n \log p - \gamma n \log n$. These
two results lead to the claimed bound on $\log p$, which concludes the proof.

\subsection{\texorpdfstring{Proof of Theorem \protect\ref{Thm4}}{Proof of Theorem 4}} \label{SecA8}
Let us fix an arbitrary subset $\{\mathbf{x}_{j_1},\ldots,
\mathbf{x}_{j_s}\}$ and
denote by $\mathcal{A}_s^{p-s}$ the set of $s$-subspaces spanned by $s$
of the remaining $p-s$ $\mathbf{x}_j$'s. By assumption, $\mathcal
{A}_s^{p-s}$
lies in a neighborhood in the Grassmann manifold $G_{n, s}$ that is
characterized by the set $R_{\delta_1} = \{(x_1,\ldots, x_s) \in[0,
1]^s\dvtx  \min_{i = 1}^s x_i \geq\delta_1^2\}$, since $x_i = \sin^2
\theta
_i$. In view of (\ref{003}), a change of variable $y_i = (1-\delta
_1^2)^{-1} (1-x_i)$ gives
%
%
\begin{eqnarray}
\label{015} d \nu&=& \bigl(1-\delta_1^2
\bigr)^{s^2/2} K_{n, s} f_1(y_1,\ldots,
y_s) \nonumber\\[-8pt]\\[-8pt]
&&{}\times\prod_{1
\leq i < j \leq s} |y_i -
y_j| \prod_{i = 1}^s
y_i^{-1/2} \,d y_1 \cdots d y_s,\nonumber
\end{eqnarray}
where $f_1(y_1,\ldots, y_s) = \prod_{i = 1}^s [1 - (1-\delta_1^2)
y_i]^{\alpha-1}$ over $[0, 1]^s$. Clearly $[1 - (1-\delta_1^2)
y_i]^{\alpha-1} \geq(1 - y_i)^{\alpha-1}$ for $\alpha\geq1$, which
together with (\ref{003}) and (\ref{015}) entails that $(1-\delta
_1^2)^{s^2/2}$ is a lower bound on the integral $\nu(R_{\delta_1})$.
However, we need an upper bound on it. The idea is to bound the
function $f_1$ by an exponential function.

We are more interested in the asymptotic behavior of $\nu(R_{\delta
_1})$ when $\delta_1$ is near zero. Using the inequality $\log(1+t)
\leq t$, we derive
\[
1 - \bigl(1-\delta_1^2\bigr) y = \bigl(1-
\delta_1^2\bigr)\bigl[1 + \bigl(1-\delta_1^2
\bigr)^{-1} \delta_1^2 -y\bigr] \leq\bigl(1-
\delta_1^2\bigr) e^{(1-\delta_1^2)^{-1} \delta_1^2} e^{-y}.
\]
This leads to
\[
f_1(y_1,\ldots, y_s) \leq\bigl[\bigl(1-
\delta_1^2\bigr) e^{(1-\delta_1^2)^{-1}
\delta
_1^2}\bigr]^{s(\alpha-1)}
\prod_{i = 1}^s e^{-(\alpha-1) y_i}.
\]
Thus we have
%
%
\begin{equation}
\label{016} \nu(R_{\delta_1}) \leq\bigl(1-\delta_1^2
\bigr)^{2^{-1}s^2+s(\alpha-1)} e^{(1-\delta_1^2)^{-1} \delta_1^2
s(\alpha-1)} K_{n, s} B_\alpha,
\end{equation}
where $B_\alpha= \int_{[0, 1]^s} \prod_{1 \leq i < j \leq s} |y_i -
y_j| \prod_{i = 1}^s y_i^{-1/2} e^{-(\alpha-1) y_i} \,d y_1 \cdots d
y_s$. Since $\alpha-1= (n-2s-1)/2 \rightarrow\infty$ as $n
\rightarrow
\infty$, a change of variable $z_i = (\alpha-1) y_i$ gives
\begin{eqnarray*}
B_\alpha& = &(\alpha-1)^{-s^2/2} \int_{[0, \alpha-1]^s}
\prod_{1
\leq i
< j \leq s} |z_i - z_j|
\prod_{i = 1}^s z_i^{-1/2}
e^{-z_i} \,d z_1 \cdots d z_s
\\
& \leq &(\alpha-1)^{-s^2/2} \int_{[0, \infty)^s} \prod
_{1 \leq i < j
\leq s} |z_i - z_j| \prod
_{i = 1}^s z_i^{-1/2}
e^{-z_i} \,d z_1 \cdots d z_s.
\end{eqnarray*}
Note that the last integral is a Selberg type integral related to the
Laguerre polynomials [\citet{Meh04}], which can be calculated exactly.
This along with the identities $\Gamma(t + 1)=t\Gamma(t)$ and $\Gamma
(3/2)=\pi^{1/2}/2$ yields
\[
B_\alpha\lesssim(\alpha-1)^{-s^2/2} s! \pi^{-s/2} \prod
_{i = 1}^s \Gamma^2(i/2).
\]

By assumption, $s \sim\gamma n$ as $n \rightarrow\infty$. It is easy
to show that
%
%
\begin{equation}
\label{017} \log\bigl[\bigl(1-\delta_1^2
\bigr)^{2^{-1}s^2+s(\alpha-1)} e^{(1-\delta_1^2)^{-1}
\delta_1^2 s(\alpha-1)}\bigr] \sim-c_{\delta_1} \gamma
n^2 + O(n),
\end{equation}
where $c_{\delta_1} = 2^{-1} [\log(1-\delta_1^2)^{-1} ]
(1-\gamma) - 2^{-1}(1-\delta_1^2)^{-1} \delta_1^2(1-2\gamma)$. It
remains to consider the term $K_{n, s} B_\alpha$. By (\ref{002}), we have
\begin{eqnarray*}
&& 2^{-s} \pi^{-s/2} \widetilde{K}_{n, s} \prod
_{i = 1}^s \Gamma^2(i/2) \\
&&\qquad= \prod
_{i=0}^{s-1}\frac{\Gamma((n-i)/2)}{\Gamma((n-s-i)/2)}
\\
&&\qquad \sim\prod_{i=0}^{s-1} \biggl(
\frac{n-i-2}{n-s-i-2} \biggr)^{(n-s-i-1)/2} \biggl[\frac
{(n-2)!}{(n-s-2)!}
\biggr]^{s/2} (2e)^{-s^2/2},
\end{eqnarray*}
where we used Stirling's formula for the Gamma function in the last
step. It follows from $s \sim\gamma n$ that
\[
\prod_{i=0}^{s-1} \biggl(\frac{n-i-2}{n-s-i-2}
\biggr)^{(n-s-i-1)/2} \lesssim\biggl(\frac{1-\gamma}{1-2\gamma} \biggr
)^{\gamma(1-\gamma)n^2/2},
\]
and an application of Stirling's formula for large factorials gives
\begin{eqnarray*}
\biggl[\frac{(n-2)!}{(n-s-2)!} \biggr]^{s/2} (2e)^{-s^2/2} & \sim &
\biggl(\frac{n-2}{n-s-2} \biggr)^{s(n-s-3/2)/2} \bigl[(n-2)/\bigl(2e^2
\bigr)\bigr]^{s^2/2}
\\
& \lesssim &(1-\gamma)^{-\gamma(1-\gamma)n^2/2} \bigl[n/\bigl(2e^2\bigr)
\bigr]^{\gamma^2 n^2/2}.
\end{eqnarray*}
Note that $(\alpha-1)^{-s^2/2} \sim[(1-2\gamma) n/2]^{-\gamma
^2n^2/2}$. Combining these results together yields
%
%
\begin{eqnarray}
\label{018} \log(K_{n, s} B_\alpha) & = & \log\Biggl[(
\alpha-1)^{-s^2/2} 2^{-s} \pi^{-s/2} \widetilde{K}_{n, s}
\prod_{i = 1}^s \Gamma^2(i/2)
\Biggr]
\\
& \lesssim & -\bigl[2\gamma+\log(1-2\gamma)\bigr] \gamma n^2/2.
\nonumber
\end{eqnarray}
It follows from (\ref{016})--(\ref{018}) that
%
%
\begin{equation}
\label{019} \log\bigl(\nu(R_{\delta_1})\bigr) \lesssim-c_{\delta_1}
\gamma n^2 -\bigl[2\gamma+\log(1-2\gamma)\bigr] \gamma
n^2/2 + O(n).
\end{equation}

Finally, we are ready to derive a bound on the dimensionality $p$.
Since $\mathcal{A}_s^{p-s}$ lies in a neighborhood in $G_{n,s}$
characterized by the set $R_{\delta_1}$ and all the subspaces in
$\mathcal{A}_s^{p-s}$ have maximum chordal distance
at least $\delta$, it holds that ${p-s \choose s} = |\mathcal
{A}_s^{p-s}| \leq\nu(R_{\delta_1})/\nu(B_{\delta, d_m})$. Aided by
(\ref{013}) and (\ref{019}), a similar argument as in the proof of
Theorem \ref{Thm3} gives the claimed bound on $\log(p-s)$. This
completes the proof.

\subsection{\texorpdfstring{Proof of Theorem \protect\ref{Thm5}}{Proof of Theorem 5}} \label{SecA9}
To prove the conclusion, we use the terminology introduced in Section~\ref{Sec3}. Note that the $n$-vectors $\mathbf{x}_j$ and $\by$ can
be viewed
as elements of Grassmannian manifold $G_{n, 1}$, which consists of all
one-dimensional subspaces of $\mathbb{R}^n$. The absolute correlation
between two $n$-vectors is given by $\cos\theta_1$, where $\theta_1
\in[0, \pi/2]$ is the principal angle between the two corresponding
one-dimensional subspaces. We use the parametrization with local
coordinate $\theta_1$ at the one-dimensional subspace $L = \{t \by\dvtx  t
\in\mathbb{R}\}$ spanned by $\by$. Then the uniform distribution on
the Grassmann manifold $G_{n, 1}$ can be expressed in local coordinate
$\theta_1$ and gives a probability measure $\nu$ in (\ref{003}) with $s
= 1$ on $[0, 1]$ through\vspace*{1pt} a change of variable $x_1 = \sin^2 \theta_1$,
where $\prod_{1 \leq i < j \leq s} |x_i - x_j| = 1$ in this case.
Consider\vspace*{2pt} the maximum chordal distance
on $G_{n, 1}$, which is defined as $\sin\theta_1$. For any $t > 0$,
denote by $B_{t, d_m}$ a ball of radius $t$ centered at $L$ in $G_{n,
1}$ under the maximum chordal distance,
that is, $B_{t, d_m} = \{x_1 \in[0, 1]\dvtx  x_1 = \sin^2 \theta_1 \leq t^2\}$
in local coordinate. We need to calculate the volumes of $B_{t_1, d_m}$
with $t_1 = 2^{-1} \sin\cos^{-1} \delta= (1-\delta^2)^{1/2}/2$ and
$B_{t_2, d_m}^c$, the complement of $B_{t_2, d_m}$ with $t_2 = \sin
\cos
^{-1} r = (1-r^2)^{1/2}$, under the measure $\nu$.

In view of (\ref{003}), we have
%
%
\begin{equation}
\label{024} \nu(B_{t, d_m}) = \nu\bigl(\bigl[0, t^2\bigr]
\bigr) = K_{n, 1} \int_0^{t^2}
x_1^{(n-3)/2} (1-x_1)^{-1/2}
\,dx_1,
\end{equation}
where $K_{n, 1} = \widetilde{K}_{n, 1}/2 = A_1^2 A_{n-1}/(4 A_n)$ with
$A_j = 2 \pi^{j/2}/\Gamma(j/2)$ the area of the unit sphere $S^{j -1}
\subset\mathbb{R}^j$. It follows from Stirling's formula for the Gamma
function that
%
%
\begin{eqnarray}
\label{025} K_{n, 1} & = & \pi^{-1/2} \Gamma(n/2)/\Gamma
\bigl((n-1)/2\bigr) \nonumber\\
&\sim&\pi^{-1/2} \biggl(\frac{n-2}{2e}
\biggr)^{1/2} \biggl(\frac{n-2}{n-3} \biggr)^{(n-2)/2}
\\
& \sim & (2 \pi)^{-1/2} n^{1/2}.
\nonumber
\end{eqnarray}
It remains to evaluate the integral in (\ref{024}). Note that
$(1-x_1)^{-1/2}$ is bounded between $1$ and $\infty$ on $[0, t^2]$ for
$t$ bounded away from $1$. Thus, we have
%
%
\begin{equation}
\label{026}\quad \int_0^{t^2} x_1^{(n-3)/2}
(1-x_1)^{-1/2} \,dx_1 > \int_0^{t^2}
x_1^{(n-3)/2} \,dx_1 = 2 (n-1)^{-1}
t^{n-1},
\end{equation}
where\vspace*{1pt} both sides have the same asymptotic order. Combining (\ref
{024})--(\ref{026}) yields $\nu(B_{t, d_m}) > c_n t^{n - 1}$ with $c_n
\sim(2/\pi)^{1/2} n^{-1/2}$, and $\nu(B_{t, d_m}) \sim\break(2/\pi)^{1/2}
n^{-1/2} t^{n - 1}$ for $t$ bounded away from $1$. Since all the $p -
s$ noise predictors $\mathbf{x}_j$ have absolute correlations
bounded by
$\delta\in(0, 1)$, we have
%
%
\begin{equation}
\label{027} p - s \leq\nu\bigl(B_{t_2, d_m}^c\bigr)/
\nu(B_{t_1, d_m}),
\end{equation}
if there exists no noise predictor that has absolute correlation with
$\by$ larger than $r \in(0, 1)$. The right-hand side of (\ref{027}) is
$[1 - \nu(B_{t_2, d_m})]/\nu(B_{t_1, d_m})$, which is less than and has
the same asymptotic order as $1/\nu(B_{t_1, d_m}) < c_n^{-1}
[4/(1-\delta^2)]^{(n-1)/2} \sim(\pi/2)^{1/2} n^{1/2} [4/(1-\delta
^2)]^{(n-1)/2}$. This together with (\ref{027}) concludes the proof.

\section{Geometry and invariant measure of Grassmann manifold} \label{Sec31}
We briefly introduce some necessary background and terminology on the
geometry and invariant measure of Grassmann manifold. Let $V_1$ and
$V_2$ be two $s$-dimensional subspaces of $\mathbb{R}^n$ and
$d_g(\cdot,\cdot) = \arccos|\langle\cdot, \cdot\rangle|$ be the geodesic
distance on $S^{n - 1}$, that is, the distance induced by the Euclidean
metric on $\mathbb{R}^n$. It was shown by \citet{Jam54} that as $v_1$
and $v_2$ vary over $V_1 \cap S^{n - 1}$ and $V_2 \cap S^{n - 1}$,
respectively, $d_g(v_1, v_2)$ has a set of $s$ critical values $\angle
(V_1, V_2) = (\theta_1,\ldots, \theta_s)$ with $\pi/2 \geq\theta_1
\geq\cdots\geq\theta_s \geq0$, corresponding to $s$ pairs of unit
vectors $(v_{1i}, v_{2i})$, $i = 1,\ldots, s$. Each critical value
$\theta_i$ is exactly the angle between $v_{1i}$ and $v_{2i}$, and
$v_{1i}$ is orthogonal to $v_{1j}$ and $v_{2j}$ if $j \neq i$. The
principal angles $\theta_i$ are unique and if none of them are equal,
the principle vectors $(v_{1i}, v_{2i})$ are unique up to a
simultaneous direction reversal. In general, the dimensions of $V_1$
and $V_2$ can be different, in which case $s$ should be their minimum.

All $s$-dimensional subspaces of $\mathbb{R}^n$ form a space, the
so-called Grassmann manifold $G_{n, s}$. It is a compact Riemannian
homogeneous space, of dimension $s (n - s)$, isomorphic to $O(n)/(O(s)
\times O(n - s))$, where $O(j)$ denotes the orthogonal group of order
$j$. It is well known that $G_{n, s}$ admits an invariant measure $\mu
$. It can be constructed by viewing $G_{n, s}$ as $V_{n, s}/O(s)$,
where $V_{n, s} \cong O(n)/O(n - s)$ denotes the Stiefel manifold of
all orthonormal $s$-frames (i.e., sets of $s$ orthonormal vectors)
in~$\mathbb{R}^n$. By deriving the exterior differential forms on those
manifolds [\citet{Jam54}], $d \mu(V)$ can be expressed in local
coordinates, at the $s$-dimensional subspace with generator matrix
$(I_s\enskip 0)$, as a product of three independent densities $\prod_{i =
1}^3 \,d \mu_i$, where
%
%
\begin{equation}
\label{001} d \mu_1 = \widetilde{K}_{n, s} \prod
_{i = 1}^s (\sin\theta_i)^{n - 2
s}
\prod_{1 \leq i < j \leq s} \bigl(\sin^2
\theta_i - \sin^2 \theta_j\bigr) \,d
\theta_1 \cdots d \theta_s
\end{equation}
over $\Theta= \{(\theta_1,\ldots, \theta_s)\dvtx  \pi/2 > \theta_1 >
\cdots> \theta_s > 0\}$, and $d \mu_2$ and $d \mu_3$ are independent
of parameters $(\theta_1,\ldots, \theta_s)$.
The normalization constant is given by
%
%
\begin{equation}
\label{002} \widetilde{K}_{n, s} = \prod_{i = 0}^{s - 1}
\frac{A_{s - i}^2 A_{n
- s
- i}}{2 A_{n - i}},
\end{equation}
where $A_j = 2 \pi^{j/2}/\Gamma(j/2)$ is the area of the unit sphere
$S^{j -1}$. A change of variable $x_i = \sin^2 \theta_i$ and
symmetrization in (\ref{001}) yield a probability measure $\nu$ on $[0,
1]^s$ with density
%
%
\begin{equation}
\label{003} d \nu= K_{n, s} \prod_{1 \leq i < j \leq s}
|x_i - x_j| \prod_{i =
1}^s
x_i^{\alpha- 1} \prod_{i = 1}^s
(1 - x_i)^{-1/2} \,d x_1 \cdots d
x_s,
\end{equation}
where $K_{n, s} = \widetilde{K}_{n, s}/(2^s s!)$ and $\alpha= (n - 2 s
+ 1)/2$.
\end{appendix}

\section*{Acknowledgments}

The author sincerely thanks the Co-Editor, Associate Editor and two
referees for their valuable comments that improved significantly the
paper.



\printaddresses

\end{document}